\documentclass[reqno,letterpaper,draft,11pt]{amsart}

\newcommand{\abs}[1]{\lvert #1 \rvert}
\begin{document}
\title{Contact  angle for immersed  surfaces in $S^{2n+1}$}
\author{Rodrigo Ristow Montes}
\address{}
\curraddr{Instituto de Matem\'atica e Estat\'{\i}stica da Universidade
  de S\~ao Paulo, Caixa Postal 66281, CEP 05315-970, S\~ao Paulo, SP,
  Brazil.}
\email{rrm@ime.usp.br}
\author{Jose  A. Verderesi}
\address
{Instituto de Matem\'atica e Estat\'{\i}stica - Universidade
  de S\~ao Paulo,
   Caixa Postal 66281, CEP 05315-970, S\~ao Paulo , SP,
  Brazil.}
\email{javerd@ime.usp.br}

\date{}
\begin{abstract}
\noindent
  In this paper  we introduce the notion of contact angle.  We deduce formulas
  for  Laplacian and  Gaussian curvature of  a minimal surface  in $S^{2n+1}$
  and  give a characterization of  the generalized Clifford Torus  as the only
  non-legendrian minimal surface in $S^5$ with constant Contact  and  K$\ddot{\mbox{a}}$hler angles.
\end{abstract}
\maketitle


\section{Introduction}

 The  K$\ddot{\mbox{a}}$hler angle was introduced by  Chern and Wolfson
 in \cite{CW} and \cite{W} and it is the most fundamental invariant 
for  minimal surfaces in complex manifolds. Using the technique of
 moving frames Wolfson obtained equations for the Laplacian and Gaussian
 curvature for an immersed minimal surface in $\mathbb{CP}^n$.

 Later, Kenmotsu in \cite{K}, Ohnita in \cite{O} and  Ogata in
 \cite{Og} classified  minimal surfaces with constant Gaussian curvature and
 constant  K$\ddot{\mbox{a}}$hler angle.

 A few years ago, Zhenqui in  \cite{Z} gave a counterexample for the Bolton's conjecture, in  \cite{BJRW}, that says that a minimal
 immersion (non-holomorphic, non anti-holomorphic, non totally real) of
 two-sphere in $\mathbb{CP}^n$ with constant  K$\ddot{\mbox{a}}$hler angle have constant Gaussian curvature.

In  \cite{MV} we introduced the contact angle and it has can be consider   a new geometric invariant in order to investigate immersed
surfaces in $S^{3}$. Geometrically, the contact angle is the complementary
angle between the contact distribution and the tangent space of the
surface. Also in  \cite{MV}, we deduced Gaussian curvature
and Laplacian formulae for an immersed minimal surface in $S^3$ and gave a
characterization of the  Clifford Torus as the only minimal surface in $S^3$ with constant Contact angle.

In this work, a  corresponding  statement to Bolton's conjecture (see
 \cite{BJRW}) for an odd dimmensional sphere has been proved. We deduce Gaussian
curvature and Laplacian formulae for an immersed  minimal surface in $S^{2n+1}$:
\begin{eqnarray}
K & = &  -\sec^{2}\beta\abs{\nabla\beta}^2-\tan\beta\Delta\beta
-2\cos\alpha(1+2\tan^{2}\beta)\beta_1\nonumber\\
    &&  +2\tan\beta\sin\alpha\alpha_1-4\tan^{2}\beta\cos^{2}\alpha\nonumber
\end{eqnarray}
\begin{eqnarray}
\tan\beta\Delta\beta &  = &
              -1-\tan^2\beta(\abs{\nabla\beta}^2+4\cos\alpha\beta_1)+2\tan\beta\sin\alpha\alpha_1\nonumber \\
                        &&
                        -\cos^2\alpha(4\tan^2\beta-1)+(1+\sin^2\beta)\abs{\frac{d\alpha}{2}+ w_1^2}^2\nonumber\\
                        &&  + \sum_{\lambda =3}^{2n}
   (h_{11}^{\lambda}h_{22}^{\lambda}-(h_{12}^{\lambda})^2)\nonumber
\end{eqnarray}
Using  Gaussian curvature, we characterize the Clifford Torus in $S^{5}$:\\

{\bf Theorem 1} \qquad  The Clifford Torus is the only non-legendrian  minimal
surface in $S^{5}$ with constant Contact and K$\ddot{\mbox{a}}$hler
angles.\\

When the  K$\ddot{\mbox{a}}$hler angle is null, we have an interesting
characterization of the Clifford Torus without suppose that the Contact angle is constant:\\

{\bf Theorem 2} \qquad  The Clifford Torus is the only  non-legendrian
minimal surface in $S^{5}$ with  Contact angle  $ 0 \leq  \beta <
\frac{\pi}{2} $ and null K$\ddot{\mbox{a}}$hler angle.
\section{Preliminaries}

Consider in $C^{n+1}$ the Hermitian product:
\begin{eqnarray}
(z,w)=\sum_{j=0}^n z^j\bar{w}^j\nonumber
\end{eqnarray}
The inner  product is given by:
\begin{eqnarray}
\langle z,w \rangle = Re (z,w)\nonumber
\end{eqnarray}
The unitary sphere
\begin{eqnarray}
S^{2n+1}=\{z\in\mathbb{C}^{n+1} | (z,z)=1\}\nonumber
\end{eqnarray}
The Reeb field in $S^{2n+1}$ is given by:
\begin{eqnarray}
\xi(z)=iz\nonumber
\end{eqnarray}
The usual contact distribution in  $S^{2n+1}$ is orthogonal to $\xi$:
\begin{eqnarray}
\Delta_z=\{v\in T_zS^{2n+1} | \langle \xi , v \rangle = 0 \}\nonumber
\end{eqnarray}
$\Delta$ is invariant by the complex structure of $\mathbb{C}^{n+1}$.\\
A unitary frame $(f_0, f_1, ..., f_n)$ of $\mathbb{C}^{n+1}$  induces an adapted
frame of $S^{2n+1}$, where $f_0=z$ is the position vector in $S^{2n+1}$ and 
$(f_1, ..., f_n)$ is the complex basis of $\Delta$.\\
Structure equations of $\mathbb {U}(n+1)$, $df_j=\psi_j^kf_k$  where
$\psi_j^k+\bar{\psi}_k^j=0$, can be interpreted as a mobile frame in $S^{2n+1}$
\begin{eqnarray}
dz = \psi_0^1f_1+ . . .  +\psi_0^nf_n + \theta \xi\nonumber
\end{eqnarray}
where $\theta=-i\psi_0^0$ is a real form.\\
$(\psi_0^j,\theta)$ constitute a coframe in $S^{2n+1}$ and unitary in
$\Delta$.\\
Taking:
\begin{equation}
\begin{array}{lcl}\label{eq:estru1}
\psi_0^j  &  =  &    w^j +i w^{j+n}\\
w^{2n+1}  &  =  &    \theta 
\end{array}
\end{equation}
we obtain a coframe $(w^{\lambda})_{\lambda=1, ... , 2n+1}$.\\
($\psi_j^k$)  forms satisfy structure equations of $\mathbb{U}(n)$:
\begin{equation}
\begin{array}{lclcl}\label{eq:estru2}
d\psi_j^k & + &  \psi_r^k \wedge  \psi_j^r  &  =  &  0\\
\psi_j^k  & + &  \bar{\psi_k^j}             &  =  &  0
\end{array}
\end{equation}
Taking $j=0$ in (\ref{eq:estru1}), we obtain that:
\begin{equation}
\begin{array}{lclclclcl}
dw^j & + & w_k^j \wedge w^k  &  - & w_k^{j+n} \wedge w^{k+n} &  - &  w^{j+n} \wedge w^{2n+1} & = & 0 \nonumber\\
dw^{j+n} & + & w_k^{j+n} \wedge w^k &  + &  w_k^{j} \wedge w^{k+n} &  + &  w^{j+n} \wedge
w^{2n+1} &  = &  0 \nonumber\\
dw^{2n+1} &  - &  2w^{j} \wedge w^{j+n} &  = & 0\nonumber
\end{array}
\end{equation}
Therefore Riemannian conexion forms satisfy:
\begin{equation}
\begin{array}{lcl}\label{eq:estru3}
w_{k+n}^{j+n} = \phantom{-}  w_k^j \\
w_{k+n}^{j}   = - w_k^{j+n}\\
w_{2n+1}^{j}     =  - w^{j+n} \\
w_{2n+1}^{j+n}   =  \phantom{-}  w^{j}
\end{array}
\end{equation}
For $j=1, ... , n$ the equation ( \ref {eq:estru2} ) gives second structure
equations:
\begin{equation}
\begin{array}{lclcl}
dw_k^j     &  + & \sum_{r=1}^{2n+1} w_r^j    \wedge w_k^r  & =   & w^j \wedge w^k \nonumber\\
dw_k^{j+n} &   + &  \sum_{r=1}^{2n+1} w_r^{j+n}\wedge w_k^r  &  =  &  w^{j+n} \wedge w^k\nonumber\\
dw_j^{2n+1} &    + & \sum_{r=1}^{2n+1} w_r^{2n+1}\wedge w_j^r &   = &   w^{2n+1} \wedge w^j\nonumber
\end{array}
\end{equation}
Therefore we conclude that curvature forms of $S^{2n+1}$ are:
\begin{equation}
\begin{array}{lcl}
\Omega_j^i = w^i  \wedge w^j\nonumber
\end{array}
\end{equation}
Note that in this adapted frame conexion forms of maximum degree $w_j^{2n+1}$
and $w_{j+n}^{2n+1}$ are:
\begin{equation}
\begin{array}{lcl}
w_{j}^{2n+1}     =  \phantom{-} w^{j+n}\nonumber \\
w_{j+n}^{2n+1}   =   - w^{j}\nonumber
\end{array}
\end{equation}

\section{Contact Angle for Immersed Surfaces in  $S^{2n+1}$}

Consider $S$ an immersed surface in $S^{2n+1}$.\\
The $\bold{contact}$ $\bf{angle}$ $\beta$ is the complementary angle
between the contact distribution $\Delta$ and the tangent space $TS$ of the surface.\\
Consider $(e_1,e_2)$ a local frame of $TS$, where $e_1\in TS\cap\Delta$. Then:
\begin{eqnarray}
\cos \beta = \langle \xi , e_2 \rangle \nonumber 
\end{eqnarray}
Let $v$ the projection of unitary vector field $e_2$ in $\Delta$:
\begin{eqnarray}
e_2 = \sin\beta v + \cos\beta \xi  \nonumber
\end{eqnarray}
We define $\alpha$ as the angle given by:
\begin{eqnarray}
\cos \alpha = \langle ie_1 , v \rangle \nonumber
\end{eqnarray}
The angle  $\alpha$ is the corresponding to the K$\ddot{\mbox{a}}$hler angle and it was introduced by Chern and Wolfson in \cite{CW}.\\
Consider $(f_j)$ an adapted frame to the distribution $\Delta$.\\
The $(f_j)$ is given by:
\begin{equation}
\begin{array}{lcl}
f_1      &   =   &  \frac{e_1-i\upsilon}{2\cos\frac{\alpha}{2}}\nonumber\\       
f_2      &   =   &  \frac{e_1+i\upsilon}{2\sin\frac{\alpha}{2}}\nonumber\\
f_j      &  \in  &  \Delta \nonumber\\
f_{n+j}  &   =   &  i f_j , \quad j=1, ...,n \nonumber\\
f_{2n+1} &   =   &  \xi \nonumber 
\end{array}
\end{equation}
Therefore $(e_1,e_2)$ can be written, as:
\begin{eqnarray}
e_1 & = &  \cos\frac{\alpha}{2}f_1+\sin\frac{\alpha}{2}f_2\nonumber\\ 
e_2 & = &  \sin\beta(\cos\frac{\alpha}{2}f_{n+1}-\sin\frac{\alpha}{2}f_{n+2})+ \cos\beta \xi\nonumber 
\end{eqnarray}
Consider normal vector fields:
\begin{equation}\label{eq:normals}
\begin{array}{lcl}
e_{n+1}     & = & \phantom{-}\sin\frac{\alpha}{2}f_1-\cos\frac{\alpha}{2}f_2\\
e_{n+2}     & = & \phantom{-}\sin\frac{\alpha}{2}f_{n+1}+\cos\frac{\alpha}{2}f_{n+2}\\
e_{2n+1}    & = &  -\cos\beta(\cos\frac{\alpha}{2}f_{n+1}-\sin\frac{\alpha}{2}
f_{n+2})+\sin\beta\xi\\
e_\lambda   & = &  \phantom{-}f_\lambda; \quad  3 \leq \lambda \leq {n} \quad  , \quad {n+3} \leq \lambda \leq {2n}
\end{array}
\end{equation}
where $(e_{\lambda})_{\lambda=1, ..., 2n+1}$ constitute a Darboux frame of $S$.\\
Let $(w^j)$, $(\theta^j)$ respective coframes of $(f_j)$ and $(e_j)$.\\
When we restrict to $S$, we obtain:
\begin{equation} \label{eq:restrict}
\begin{array}{lcl}
w^{1}       & = &  \phantom{-}\cos\frac{\alpha}{2}\theta^1\\
w^{2}       & = &  \phantom{-}\sin\frac{\alpha}{2}\theta^1\\ 
w^{n+1}     & = &  \phantom{-}\cos\frac{\alpha}{2}\sin\beta\theta^2\\ 
w^{n+2}     & = &  -\sin\frac{\alpha}{2}\sin\beta\theta^2\\
w^{2n+1}    & = &  \phantom{-}\cos\beta\theta^2\\ 
w^{\lambda} & = &  \phantom{-}0  \hspace{2cm}  3 \leq \lambda \leq {n} \quad  , \quad {n+3} \leq \lambda \leq {2n}
\end{array}
\end{equation}
\section{Equations for curvature and Laplacian of a minimal surface in $S^{2n+1}$} \mbox{}
Consider $D$ the covariant derivative of Riemannian metric in
$S^{2n+1}$. Then:
\begin{equation}
\begin{array}{ccc}
Df_{j}  &  =  &  w_j^k f_k \nonumber
\end{array}  
\end{equation}
where $(w_j^k)$ satisfy (\ref{eq:estru3}).\\ 
Consider conexion forms $(\theta_k^j)$ associate to the Darboux frame
$(e_j)$. Then:
\begin{equation}
\begin{array}{ccc}
De_{j} = \theta_j^k e_k \nonumber
\end{array}  
\end{equation}
Second fundamental forms in this frame are given by:
\begin{equation}
\begin{array}{lclll}
II^{n+1}      & = & \theta_1^{n+1} \theta^1 &   +   &  \theta_2^{n+1} \theta^2\\
II^{n+2}      & = & \theta_1^{n+2} \theta^1 &   +  &  \theta_2^{n+2} \theta^2\\
II^{2n+1}     & = & \theta_1^{2n+1} \theta^1 &  + & \theta_2^{2n+1} \theta^2\\
II^{\lambda}  & = & \theta_1^{\lambda} \theta^1 &   + & \theta_2^{\lambda} \theta^2\\
\end{array}  
\end{equation}
Using  formulae (\ref{eq:normals}) and (\ref{eq:restrict}), we deduce:
\begin{equation}
\begin{array}{lcl}
\theta_{n+1}^1   &  =  & \phantom{-}\frac{d\alpha}{2}+w_1^2 \nonumber\\
\theta_{n+1}^2   &  =  & \phantom{-}\sin\beta(\frac{\sin\alpha}{2}w_1^{n+1}-\sin^{2}\frac{\alpha}{2}w_1^{n+2}-\cos^{2}\frac{\alpha}{2}w_2^{n+1}+\frac{\sin\alpha}{2}w_2^{n+2})\nonumber\\
\theta_{n+2}^1   &  =  & \phantom{-}\frac{\sin\alpha}{2}w_{n+1}^1+\sin^{2}\frac{\alpha}{2}
w_{n+1}^2+\cos^{2}\frac{\alpha}{2}w_{n+2}^1+\frac{\sin\alpha}{2}w_{n+2}^2 \nonumber\\
\theta_{n+2}^2   &  =  & \phantom{-}\sin\beta(\frac{d\alpha}{2}-w_{n+1}^{n+2})\nonumber\\
\theta_{2n+1}^1   &  =  & -\cos\beta(\cos^{2}\frac{\alpha}{2}w_{n+1}^1+\frac{\sin\alpha}{2}w_{n+1}^2-\frac{\sin\alpha}{2}w_{n+2}^1-\sin^{2}\frac{\alpha}{2}
w_{n+2}^2)\nonumber\\
                  &     & -\cos\alpha\sin^{2}(\beta)\theta^2 \nonumber \\
\theta_{2n+1}^2   &  =  & \phantom{-}d\beta+\cos(\alpha)\theta^{1}\nonumber
\end{array}
\end{equation}
Minimallity and symmetry conditions are:
\begin{equation}\label{eq:minimal}
\begin{array}{lcl}
\theta_1^{\lambda} \wedge \theta^1 + \theta_2^{\lambda} \wedge \theta^2 &  = & 0  \\
\theta_1^{\lambda} \wedge \theta^1 - \theta_2^{\lambda} \wedge \theta^2 &  = &0  
\end{array}
\end{equation}
From (\ref{eq:minimal}), we obtain:
\begin{equation}\label{eq: segunda}
\begin{array}{ccl}
\theta_{n+1}^1   &  =  &\phantom{-}\frac{d\alpha}{2}+w_1^2\\
\theta_{n+1}^2   &  =  &\phantom{-}\frac{d\alpha\circ J}{2}+w_1^2\circ J \\
\theta_{n+2}^1   &  =  &\phantom{-}\sin\beta(-\frac{d\alpha\circ J}{2}-w_{1}^{2}\circ
J)\\
\theta_{n+2}^2   &  =  &\phantom{-}\sin\beta(\frac{d\alpha}{2}+w_{1}^{2})\\
\theta_{2n+1}^1  &  =  & -d\beta\circ J+\cos\alpha\theta^2 \\
\theta_{2n+1}^2  &  =  &\phantom{-}d\beta+\cos\alpha\theta^{1}
\end{array}   
\end{equation} 
where the complex structure $J$  of $S$ is given by $Je_1=e_2$ and $Je_2=-e_1$.\\
Gauss equation is:
\begin{equation}\label{eq: Gauss}
\begin{array}{lcl}
d\theta_2^1 + \theta_k^1 \wedge \theta_2^k & = & \theta^1 \wedge \theta^2
\end{array}  
\end{equation}
We define:
\begin{equation}\label{eq:soma}
\begin{array}{lcl}
\theta_j^{\lambda} & = & h_{jk}^{\lambda}\theta^k
\end{array}
\end {equation}
Then using (\ref{eq: segunda}), (\ref{eq: Gauss}) and (\ref{eq:soma}), we obtain that:
\begin{eqnarray}\label{eq:curvatura1}
K  & = &     1- \abs{\nabla(\beta)}^2-2\beta_1\cos\alpha-\cos^2\alpha\nonumber\\
   &   &
   -(1+\sin^2\beta)(\frac{\abs{\nabla\alpha}^2}{4}+w_1^2(e_1)\alpha_1+w_1^2(e_2)\alpha_2+\abs{\nabla
   w_1^2}^2)\nonumber\\
   &   & + \sum_{\lambda =3}^{2n}(h_{11}^{\lambda}h_{22}^{\lambda}-(h_{12}^{\lambda})^2)
\end{eqnarray}
Or, to be equivalent to:
\begin{eqnarray}\label{eq:kurvatura}
K   =  &&   1- \abs{\nabla\beta+\cos\alpha
  e_1}^2-(1+\sin^2\beta)\abs{\frac{d\alpha}{2}+ w_1^2}^2\nonumber\\
       &&   +     \sum_{\lambda =3}^{2n}(h_{11}^{\lambda}h_{22}^{\lambda}-(h_{12}^{\lambda})^2)
\end{eqnarray}
Taking the differential of $(\theta^1,\theta^2)$, we have:
\begin{eqnarray}
d\theta^{1} & - & \sin\beta[\cos^{2}\tfrac{\alpha}{2}w_1^{n+1} + 
            \tfrac{\sin\alpha}{2}(w_2^{n+1}  - w_1^{n+2})\nonumber\\
            & - & \sin^{2}\tfrac{\alpha}{2} w_2^{n+2}]\wedge \theta^2   =  0\nonumber\\
d\theta^{2} & + & \sin\beta[\cos^{2}\tfrac{\alpha}{2}w_1^{n+1}+\tfrac{\sin\alpha}{2}(w_2^{n+1}-w_1^{n+2})-\sin^{2}\tfrac{\alpha}{2}w_2^{n+2}\nonumber\\
            & - & \cos\beta\cos\alpha \theta^{2}]\wedge \theta^1  = 0 \nonumber
\end{eqnarray}
Therefore:
\begin{eqnarray}
\theta_2^1 = -\sin\beta[\cos^{2}\tfrac{\alpha}{2}w_1^{n+1}+\tfrac{\sin\alpha}{2}(w_2^{n+1}-w_1^{n+2})-\sin^{2}\tfrac{\alpha}{2}
w_2^{n+2}-\cos\beta\cos\alpha \theta^{2}]\nonumber
\end{eqnarray}
Using the symmetry and minimalitty and the complex structure of $S$, we obtain that:
\begin{equation} \label{eq:conex}
\begin{array}{lcl} 
\theta_2^1  &  =  &  \tan\beta(d\beta\circ J-2\cos\alpha\theta^2)
\end{array}
\end{equation}
Computing the differential of (\ref{eq:conex}), we can get:
\begin{eqnarray}
d\theta_2^1  = && (-\sec^{2}\beta\abs{\nabla\beta}^2-\tan(\beta)\Delta(\beta)-2\cos\alpha\beta_1(1+2\tan^{2}\beta)\nonumber\\
              &&
              +2\tan\beta\sin\alpha\alpha_1-4\tan^{2}\beta\cos^{2}\alpha)\theta^1
               \wedge \theta^2
\end{eqnarray}
where $\Delta\beta \quad  = \quad tr \nabla d\beta$.\\
Therefore the Gaussian curvature is:
\begin{equation}
\begin{array}{ccl} \label{eq:curvatura2}
K = && -\sec^{2}\beta\abs{\nabla\beta}^2-\tan\beta\Delta(\beta)-2\cos\alpha\beta_1(1+2\tan^{2}\beta)+\\
    && +2\tan\beta\sin\alpha\alpha_1-4\tan^{2}\beta\cos^{2}\alpha 
\end{array}
\end{equation}
Jointing (\ref{eq:curvatura1}) and (\ref{eq:curvatura2}), we obtain a new
Laplacian equation:
\begin{eqnarray}\label{eq:lapla}
\tan\beta\Delta(\beta)= &&-1-\tan^2\beta(\abs{\nabla\beta}^2+4\cos\alpha\beta_1)+2\tan\beta\sin\alpha\alpha_1\nonumber\\
                        &&  -\cos^2\alpha(4\tan^2\beta-1)+(1+\sin^2\beta)\abs{\frac{d\alpha}{2}+w_1^2}^2 \nonumber\\
                        && + \sum_{\lambda =3}^{2n}(h_{11}^{\lambda}h_{22}^{\lambda}-(h_{12}^{\lambda})^2)
\end{eqnarray}
\section{Minimal Surfaces in $S^5$} 
In this section, we compute Contact angle for examples of minimal torus given
by Kenmotsu, in \cite{KK2}.
\subsection{Contact Angle  of Legendrian Minimal Torus}\mbox{} \\ \\
Consider the torus in $S^5$ defined by:
\begin{eqnarray*}
T^3=\lbrace(z_1,z_2,z_3)\in{C^3} | z_1\bar{z_1}=\frac{1}{3},z_2\bar{z_2}=\frac{1}{3},z_3\bar{z_3}=\frac{1}{3}\rbrace
\end{eqnarray*}
We consider the immersion:
\begin{eqnarray*}
f(u_1,u_2)=\frac{\sqrt 3}{3}(e^{iu_1},e^{iu_2},e^{-i(u_1+u_2)})
\end{eqnarray*}
Then the line vector fields are:
\begin{eqnarray}
e_1 & = &  \frac{\sqrt 2}{2}i(e^{iu_1},0,-e^{-i(u_1+u_2)})\nonumber\\ 
e_2 & = &  \frac{\sqrt 2}{\sqrt 3}i(-\frac{1}{2}e^{iu_1},e^{iu_2},- \frac{1}{2}e^{-i(u_1+u_2)})  \nonumber 
\end{eqnarray}
The normal vector field is:
\begin{equation}
\begin{array}{lcl}
\xi  &  =  & i\frac{\sqrt 3}{3}(e^{iu_1},e^{iu_2},e^{-i(u_1+u_2)})\nonumber \\
\end{array}
\end{equation}
The Contact angle  is :
\begin{equation}
\begin{array}{lcl}
\cos \beta &  =  &  \langle e_2 ,\xi \rangle\nonumber \\
           &  =  &   0 \nonumber
\end{array}
\end{equation}
Therefore, 
\begin{equation}
\begin{array}{lcl}
 \beta &  =  &  \frac{\pi}{2}\nonumber
\end{array}
\end{equation}
\subsection{Contact Angle of Generalized Clifford Torus}\mbox{} \\ \\
Consider the torus in $S^5$ defined by:
\begin{eqnarray*}
T^3=\lbrace(z_1,z_2,z_3)\in{C^3} | z_1\bar{z_1}=\frac{1}{3},z_2\bar{z_2}=\frac{1}{3},z_3\bar{z_3}=\frac{1}{3}\rbrace
\end{eqnarray*}
We consider the immersion:
\begin{eqnarray*}
f(u_1,u_2)=\frac{\sqrt 3}{3}(e^{iu_1},e^{iu_2},e^{i(u_2-u_1)})
\end{eqnarray*}
Then the line vector fields are:
\begin{eqnarray}
e_1 & = &  i\frac{\sqrt 2}{2}(e^{iu_1},0,-e^{i(u_2-u_1)})\nonumber\\ 
e_2 & = &  i\frac{\sqrt 2}{\sqrt 3}(\frac{1}{2}e^{iu_1},e^{iu_2}, \frac{1}{2}e^{i(u_2-u_1)})  \nonumber 
\end{eqnarray}
The normal vector field is:
\begin{equation}
\begin{array}{lcl}
\xi  &  =  & i\frac{\sqrt 3}{3}(e^{iu_1},e^{iu_2},e^{i(u_2-u_1)})\nonumber \\
\end{array}
\end{equation}
The Contact angle  is :
\begin{equation}
\begin{array}{lcl}
\cos \beta &  =  &  \langle e_2 ,\xi \rangle\nonumber \\
           &  =  &  \frac{2\sqrt 2}{3} \nonumber
\end{array}
\end{equation}
Therefore, 
\begin{equation}
\begin{array}{lcl}
 \beta &  =  &  arc \quad cos (\frac{2\sqrt 2}{3})\nonumber
\end{array}
\end{equation}
\subsection{Contact Angle of Clifford Torus}\mbox{}\\ \\
Consider the torus in $S^5$ defined by:
\begin{eqnarray*}
T^2=\lbrace(z_1,z_2,0)\in{C^3} | z_1\bar{z_1}=\frac{1}{2},z_2\bar{z_2}=\frac{1}{2}\rbrace
\end{eqnarray*}
We consider the immersion:
\begin{eqnarray*}
f(u_1,u_2)=\frac{\sqrt 2}{2}(e^{iu_1},e^{iu_2},0)
\end{eqnarray*}
Then the line vector fields are:
\begin{eqnarray}
e_1 & = &  i\frac{\sqrt 2}{2}(e^{iu_1},-e^{iu_2},0)\nonumber\\ 
e_2 & = &  i\frac{\sqrt 2}{2}(e^{iu_1},e^{iu_2},0)  \nonumber 
\end{eqnarray}
The normal vector field is:
\begin{equation}
\begin{array}{lcl}
\xi  &  =  & i\frac{\sqrt 2}{2}(e^{iu_1},e^{iu_2},0)\nonumber \\
\end{array}
\end{equation}
The Contact angle  is :
\begin{equation}
\begin{array}{lcl}
\cos \beta &  =  &  \langle e_2 ,\xi \rangle\nonumber \\
           &  =  &   1 \nonumber
\end{array}
\end{equation}
Therefore, 
\begin{equation}
\begin{array}{lcl}
 \beta &  =  &  0\nonumber
\end{array}
\end{equation}

\section{Main Results}
\subsection{Proof of the Theorem 1}\mbox{}

Suppose that $\alpha $ and $\beta$ are constants, then equation
(\ref{eq:curvatura2} ) reduces to:
\begin{equation}\label{eq:curvaturasimp}
\begin{array}{ccl} 
K = && -4\tan^{2}\beta\cos^{2}\alpha 
\end{array}
\end{equation}
As the line field determined by $e_1$ is globally defined, we have therefore
that $S$ is parallelizable and by equation (\ref{eq:curvaturasimp})  Gaussian
curvature is constant.\\ 
Therefore, using Gauss-Bonnet theorem, we have that Gaussian curvature of  $S$ is null everywhere.\\
Using the work of Kenmotsu, in \cite{KK}, we have that $S$ is the
generalized Clifford Torus in $S^5$.\hfill{q.e.d.}

\subsection{Proof of the Theorem 2}\mbox{}

When the surface is immersed in $S^5$, equation (\ref{eq:kurvatura}) is :
\begin{eqnarray}\label{eq: curvatura3}
K  & = &     1- \abs{\nabla\beta+\cos\alpha e_1}^2 -\
   (1+\sin^2\beta)\abs{\frac{d\alpha}{2}+w_1^2}^2 
\end{eqnarray} 
In $S^5$ equation (\ref{eq:lapla}) is given by :
\begin{eqnarray}\label{eq:lapla2}
\tan\beta\Delta(\beta)= &&-1-\tan^2\beta(\abs{\nabla\beta}^2+4\cos\alpha\beta_1)+2\tan\beta\sin\alpha\alpha_1\nonumber\\
                        &&  -\cos^2\alpha(4\tan^2\beta-1)+(1+\sin^2\beta)\abs{\frac{d\alpha}{2}+ w_1^2}^2 \nonumber\\
\end{eqnarray} 
Suppose that $ \alpha=0 $, then the adapted frame $(f_j)$ is given by:
\begin{equation}
\begin{array}{lcl}
f_1         &   =   &  e_1\nonumber\\       
f_3         &   =   &  ie_1\nonumber\\
(f_2,f_4)   &  \in  &  \Delta \nonumber\\
f_5         &   =   &  \xi \nonumber 
\end{array}
\end{equation}
Therefore $(e_1,e_2)$ can be written, as:
\begin{eqnarray}
e_1 & = &  f_1\nonumber\\ 
e_2 & = &  \sin\beta f_3+ \cos\beta \xi\nonumber 
\end{eqnarray}
Normal vector fields are:
\begin{equation}
\begin{array}{lcl}
e_{3}     & = & \phantom{-}f_2\\
e_{4}     & = & \phantom{-}f_4\\
e_{5}    & = &  -\cos\beta f_3+\sin\beta\xi
\end{array}
\end{equation}
Restrict to $S$, we obtain:
\begin{equation} 
\begin{array}{lcl}
w^{1}       & = &  \theta^1\\
w^{2}       & = &  0\\ 
w^{3}       & = &  \sin\beta\theta^2\\ 
w^{4}       & = &  0\\
w^{5}       & = &  \cos\beta\theta^2\\ 
\end{array}
\end{equation}
Taking the differential of $e_1$, we obtain:
\begin{equation} 
\begin{array}{lcl}
De_1        & = &  Df_1\\
            & = &  w_1^2f_2+w_1^3f_3+w_1^4f_4+w_1^5f_5\\ 
            & = &  (\sin\beta
            w_1^3+\tfrac{\sin 2\beta}{2}\theta^2)e_2-w_1^2e_3+w_1^4e_4+(-\cos\beta
            w_1^3+\sin^2\beta\theta^2)e_5\\ 
\end{array}
\end{equation}
Therefore:
\begin{equation} \label{eq: formas1}
\begin{array}{lcl}
\theta_1^2      & = & \phantom{-} \sin\beta w_1^3 + \frac{\sin 2\beta}{2}\theta^2\\
\theta_1^3      & = &  -w_1^2\\ 
\theta_1^4      & = &  \phantom{-} w_1^4\\
\theta_1^5      & = &  -\cos\beta w_1^3+\sin^2\beta\theta^2\\ 
\end{array}
\end{equation}
Taking the differential of $e_2$, we obtain:
\begin{equation} 
\begin{array}{lcl}
De_2        & = &  \cos\beta d\beta f_3 - \sin\beta d\beta f_5 + \sin\beta Df_3
            + \cos\beta Df_5\\
            & = &  w_1^2f_2+w_1^3f_3+w_1^4f_4+w_1^5f_5\\ 
            & = &  (\sin\beta
            w_3^1-\tfrac{\sin 2\beta}{2}\theta^2)e_1-\sin\beta w_1^4 e_3 +
            \sin\beta w_1^2 e_4 - (d\beta+\theta^1)e_5\\ 
\end{array}
\end{equation} 
Therefore:
\begin{equation} \label{eq: formas2}
\begin{array}{lcl}
\theta_2^1      & = & \phantom{-} \sin\beta w_3^1 - \frac{\sin 2\beta}{2}\theta^2\\
\theta_2^3      & = &  -\sin\beta w_1^4\\ 
\theta_2^4      & = &  \phantom{-} \sin\beta w_1^2\\
\theta_2^5      & = &  -(d\beta+\theta^1)\\ 
\end{array}
\end{equation}
Equivalently to minimallity and symmetry conditions, we have:
\begin{equation}\label{eq: complexa} 
\begin{array}{lcl}
\theta_1^3 \circ J      & = & \theta_2^3\\
\theta_1^4 \circ J      & = & \theta_2^4
\end{array}
\end{equation}
Substituing (\ref{eq: formas1})  and (\ref{eq: formas2}) at the equations
(\ref{eq: complexa}), we obtain:
\begin{equation} 
\begin{array}{lcl}
w_1^2  \circ J      & = & \sin\beta w_1^4\\
w_1^4  \circ J      & = & \sin\beta w_1^2
\end{array}
\end{equation}
Therefore:
\begin{equation} 
\begin{array}{lcl}
(1-\sin^2\beta)w_1^2 \circ J   & = & 0\\
\end{array}
\end{equation}
The  surface $S$ is non-legendrian, we conclude that $w_1^2=0$. Then the  equation (\ref{eq:lapla2} ) reduces to:
\begin{equation}
\begin{array}{ccl} 
\Delta(\beta) &  =   & -\tan(\beta)\abs{\nabla\beta+2e_1}^2
\end{array}
\end{equation}
Suppose that $ 0\leq  \beta < \frac{\pi}{2}$, then using \textbf{ Hopf 's
  lemma} we obtain that  $\beta$ is constant and using equation
  (\ref{eq: curvatura3}), we obtain that K=0, and therefore $S$ is the Clifford
  torus. \hspace{9 cm} \hfill{q.e.d.}

\end{document}